\documentclass{amsart}

\usepackage{fullpage}
\usepackage{hyperref}
\usepackage[capitalize]{cleveref}
\usepackage{mathrsfs}
\usepackage[frozencache]{minted}
\usepackage{physics}
\usepackage{float}

\newcommand{\CC}{\mathbb C}
\newcommand{\RR}{\mathbb R}

\title{ForeignFunctions package for Macaulay2}
\author{Douglas A. Torrance}
\address{Department of Mathematical Sciences, Piedmont University, Demorest, Georgia 30535}
\email{dtorrance@piedmont.edu}

\subjclass[2020]{Primary 13-04, 14-04; Secondary 62J05, 65T50}

\begin{document}

\maketitle

\begin{abstract}
  We introduce the \texttt{ForeignFunctions} package for \textsc{Macaulay2}, which uses \textsc{libffi} to provide the ability to call functions from external libraries without needing to link against them at compile time.  As examples, we use the library FFTW to multiply polynomials using fast Fourier transforms, call a LAPACK function to solve a general Gauss-Markov linear model problem, and use JIT compilation to compute Fibonacci numbers.
\end{abstract}

\section{Introduction}

Traditionally, there have been two ways to call functions from external libraries in \textsc{Macaulay2} \cite{M2}, a software platform widely used by algebraic geometers and commutative algebraists.

One is to write some code wrapping the external library's functions.  This code might involve C++ if the library is to be used by the \textsc{Macaulay2} engine, and it would also very likely involve the D language used by the interpreter and the top-level \textsc{Macaulay2} language.  The library would then be statically or dynamically linked at compile time.  For example, the \texttt{roots} function for computing the zeros of a polynomial is implemented in this way using the library \textsc{MPSolve} \cite{MPSolve1,MPSolve2}.  The addition of a new library using this method requires rebuilding \textsc{Macaulay2} from source.

The other is to use the shell to communicate with an external program's command line interface, often creating several temporary files containing the program's input and output in the process. For example, the \texttt{Normaliz} package \cite{Normaliz-M2} communicates in this way with \textsc{Normaliz}, a software platform for discrete convex geometry \cite{Normaliz}.

However, using the library \textsc{libffi} \cite{libffi}, it is possible to create a ``foreign function interface'' where users may open a shared library and call functions from it at run time, without needing to link against the library at compile time or rely on the overhead of communicating with a command line interface via the shell.  Beginning with version 1.21, \textsc{Macaulay2} has been distributed with the \texttt{ForeignFunctions} package, written by author, which provides such an interface.

This paper is organized as follows.  In \cref{section: interface}, we describe the interface.  In \cref{section: implementation}, we discuss how it was implemented.  The last three sections are devoted to examples.  In \cref{section: fft example}, we call functions from the \textsc{FFTW} library \cite{FFTW} to compute fast Fourier transforms to multiply polynomials.  In \cref{section: glm example}, we use LAPACK to solve a general Gauss-Jordan linear model problem.  And in \cref{section: jit example}, we use the \texttt{ForeignFunctions} package for JIT compilation to quickly compute Fibonacci numbers.

\section{Foreign function interface}
\label{section: interface}

The \texttt{ForeignFunctions} package introduces several new types that we outline in this section.

The first of these is \texttt{SharedLibrary}.  Each \texttt{SharedLibrary} object contains a pointer to a shared library ``handle'' as returned by the C standard library function \texttt{dlopen}.  The constructor method \texttt{open\-Shared\-Library}, which takes a string representing a library name or filename as input, calls \texttt{dlopen} and returns an instance of this class.  For an example, see \cref{listing: shared library}, where we load the \textsc{FFTW} library \cite{FFTW} that will be used in \cref{section: fft example}.

\begin{listing}[h]
\begin{minted}[frame=lines,framesep=1ex]{macaulay2}
i1 : needsPackage "ForeignFunctions"

o1 = ForeignFunctions

o1 : Package

i2 : openSharedLibrary "fftw3"

o2 = fftw3

o2 : SharedLibrary
\end{minted}
\caption{Constructing a \texttt{SharedLibrary} object}
\label{listing: shared library}
\end{listing}

There are also two related abstract classes: \texttt{ForeignType} and \texttt{ForeignObject}.  \textsc{Macaulay2} users may be familiar with \texttt{PolynomialRing}/\texttt{RingElement} and \texttt{Module}/\texttt{Vector}, the abstract classes for polynomial rings/modules and their elements, respectively.  These classes work similarly.  There are a number of instances of \texttt{ForeignType} representing various C types such as \texttt{int} (32-bit signed integers), \texttt{double} (double-precision floating point numbers) and \texttt{char*} (pointers to null-terminated strings).  Each \texttt{ForeignType} object is a class, and \texttt{Foreign\-Object} is the abstract class for instances of these classes.  Each \texttt{ForeignType} object is a self-initializing class that serves as a constructor method for the appropriate \texttt{ForeignObject} class, taking a \textsc{Macaulay2} object as input and converting it to the appropriate C type.  See, for example, \cref{listing: foreign objects}.  There are a number of built-in instances of \texttt{ForeignType}, and users may create additional types for arrays, structs, and unions.

\begin{listing}[h]
\begin{minted}[frame=lines,framesep=1ex]{macaulay2}
i3 : int 5

o3 = 5

o3 : ForeignObject of type int32

i4 : double pi

o4 = 3.14159265358979

o4 : ForeignObject of type double

i5 : charstar "Hello, world!"

o5 = Hello, world!

o5 : ForeignObject of type char*
\end{minted}
\caption{Each instance of a \texttt{ForeignType} is a \texttt{ForeignObject}}
\label{listing: foreign objects}
\end{listing}

Finally, the \texttt{ForeignFunction} class, which gives the package its name, is a subclass of \texttt{FunctionClosure}.  Instances of this class convert their input from \textsc{Macaulay2} objects to the corresponding C objects (using the appropriate \texttt{ForeignType} classes), call a C function, and then return their output as a \texttt{ForeignObject}.  \texttt{ForeignFunction} objects are created using the \texttt{foreignFunction} method, which takes a \texttt{SharedLibrary} object (this is optional if the C function is available in shared library \textsc{Macaulay2} is already linked against), a string with the name of the function, the output type, and input type (or a list of input types).  See \cref{listing: puts example} for an example \texttt{ForeignFunction} object wrapping the \texttt{puts} function from the C standard library for outputting strings.  It is possible to call functions from C++ libraries as well, but their names are often mangled by the compiler, so this must be accounted for.

\begin{listing}[h]
\begin{minted}[frame=lines,framesep=1ex]{macaulay2}
i6 : puts = foreignFunction("puts", int, charstar)

o6 = puts

o6 : ForeignFunction

i7 : puts "Hello, world!"
Hello, world!

o7 = 14

o7 : ForeignObject of type int32
\end{minted}
\caption{Printing ``Hello, world!'' with \texttt{puts} as a \texttt{ForeignFunction}}
\label{listing: puts example}
\end{listing}

\section{Implementation details}
\label{section: implementation}

Several additions to the \textsc{Macaulay2} interpreter were made to support the \texttt{ForeignFunctions} package, which we describe in this section.

First, \texttt{Pointer}, a new subclass of \texttt{Thing}, was added to represent pointers to memory.  Although users of the package often do not need to deal with these objects directly, they are fundamental to the inner workings of the package.  For example, all instances of the new classes mentioned in \cref{section: interface} are basic lists or hash tables containing \texttt{Pointer} objects.

Next, a number of low-level routines for converting \textsc{Macaulay2} objects into pointers to C objects and back again were added.  In particular, \textsc{Macaulay2} integers, which belong to the class \texttt{ZZ}, use GMP, the GNU Multiple Precision Arithmetic Library \cite{gmp}.  Functions were added to convert these to and from pointers to 8-bit, 16-bit, 32-bit, and 64-bit signed and unsigned integers.  Similarly, \textsc{Macaulay2} reals, or \texttt{RR} objects, use MPFR, the GNU Multiple Precision Floating-Point Reliable Library \cite{mpfr}.  Functions were added to convert these to and from pointers to single- and double-precision floating point numbers (\texttt{float} and \texttt{double} objects, respectively).  The constructor methods for each \texttt{ForeignType} object use these functions to convert from \textsc{Macaulay2} objects to pointers to C objects, and the \texttt{value(ForeignObject)} method uses them to convert back.

Everything is then driven by the \textsc{libffi} library \cite{libffi}.  When a \texttt{ForeignFunction} object is created using the \texttt{foreignFunction} constructor method, the function is found using the C standard library function \texttt{dlsym} and an \texttt{ffi\_cif} struct is created with information about the function to be called, its input types, and its output types.  Then when it is called, its inputs are converted to pointers to C objects, and \texttt{ffi\_call} is called to actually call the function.

Much of this code may be found in the file \texttt{M2/Macaulay2/d/ffi.d} at the \textsc{Macaulay2} GitHub repository (\url{https://github.com/Macaulay2/M2}).

\section{Fast Fourier transform example}
\label{section: fft example}

Consider $\vb a = (a_0,\ldots,a_n)\in\CC^{n+1}$ and let $f=\sum_{k=0}^na_kx^k\in\CC[x]$.  Recall that the \textit{discrete Fourier transform} of $\vb a$ is $\mathscr F\{\vb a\} = (f(1), f(\omega^{-1}), \ldots, f(\omega^{-n}))$ and the \textit{inverse discrete Fourier transform} of $\vb a$ is $\mathscr F^{-1}\{\vb a\} = (f(1), f(\omega), \ldots, f(\omega^n))$, where $\omega=\exp\left(\frac{2\pi i}{n+1}\right)$ is a principal $(n+1)$th root of unity.  In particular, $\mathscr F\{\mathscr F^{-1}\{\vb a\}\} = \mathscr F^{-1}\{\mathscr F\{\vb a\}\} = (n + 1)\vb a$.  Na\"ively, computing a discrete Fourier transform or its inverse using their definitions has time complexity $O(n^2)$, but the famous divide-and-conquer \textit{fast Fourier transform (FFT)} algorithm of Cooley and Tukey \cite{CT} has time complexity of only $O(n\log n)$.

A well-known implementation of the FFT algorithm is the C library \textsc{FFTW} \cite{FFTW}.  Using the \texttt{Foreign\-Functions} package, it is possible to create wrappers for its functions in \textsc{Macaulay2}; see \cref{listing: wrapping fftw}.  Note that \textsc{FFTW} works with arrays of \texttt{fftw\_complex} objects, but each of these is just an array containing 2 \texttt{double}'s.  So some time is spent unpacking the real and imaginary parts from a list of $n$ \texttt{CC} objects (\textsc{Macaulay2}'s complex number type) into a list of $2n$ \texttt{RR} objects that can be converted into an array of \texttt{double}'s for use by \textsc{FFTW}.

\begin{listing}[h]
\begin{minted}[frame=lines,framesep=1ex]{macaulay2}
needsPackage "ForeignFunctions"

libfftw = openSharedLibrary "fftw3"
fftwPlanDft1d = foreignFunction(libfftw, "fftw_plan_dft_1d", voidstar,
    {int, voidstar, voidstar, int, uint})
fftwExecute = foreignFunction(libfftw, "fftw_execute", void, voidstar)
fftwDestroyPlan = foreignFunction(libfftw, "fftw_destroy_plan", void, voidstar)

fftHelper = (x, sign) -> (
    n := #x;
    inptr := getMemory(2 * n * size double);
    outptr := getMemory(2 * n * size double);
    p := fftwPlanDft1d(n, inptr, outptr, sign, 64);
    registerFinalizer(p, fftwDestroyPlan);
    dbls := splice apply(x, y -> (realPart numeric y, imaginaryPart numeric y));
    apply(2 * n, i -> *(value inptr + i * size double) = double dbls#i);
    fftwExecute p;
    r := ((2 * n) * double) outptr;
    apply(n, i -> value r_(2 * i) + ii * value r_(2 * i + 1)))

fastFourierTransform = method()
fastFourierTransform List := x -> fftHelper(x, -1)

inverseFastFourierTransform = method()
inverseFastFourierTransform List := x -> fftHelper(x, 1)
\end{minted}
\caption{\textsc{Macaulay2} wrappers around \textsc{FFTW} functions}
\label{listing: wrapping fftw}
\end{listing}

One application of the FFT algorithm of possible interest to \textsc{Macaulay2} users is $O(n\log n)$ polynomial multiplication.  Consider $f=\sum_{k=0}^ma_kx^k,g=\sum_{k=0}^nb_kx^k\in\CC[x]$ and let $\vb a = (a_0,\ldots,a_m,0,\ldots,0)$ and $\vb b=(b_0,\ldots,b_n,0,\ldots,0)$ be $(m+n+1)$-tuples.  Define $\vb c = \mathscr F\{\vb a\}\cdot\mathscr F\{\vb b\}$ using componentwise multiplication, and then it follows that $\frac{1}{m+n+1}\mathscr F^{-1}\{\vb c\}$ contains the coefficients of $fg$.  See \cite[Chapter 30]{CLRS} for more details and \cref{listing: fft poly mult} for a \textsc{Macaulay2} implementation using the \textsc{FFTW} wrappers defined above.

\begin{listing}[h]
\begin{minted}[frame=lines,framesep=1ex]{macaulay2}
fftMultiply = (f, g) -> (
    m := first degree f;
    n := first degree g;
    a := fastFourierTransform splice append(
	apply(m + 1, i -> coefficient(x^i, f)), n:0);
    b := fastFourierTransform splice append(
	apply(n + 1, i -> coefficient(x^i, g)), m:0);
    c := inverseFastFourierTransform apply(a, b, times);
    sum(m + n + 1, i -> c#i * x^i)/(m + n + 1))
\end{minted}
\caption{\textsc{Macaulay2} implementation of FFT polynomial multiplication}
\label{listing: fft poly mult}
\end{listing}

The \textsc{Macaulay2} engine multiplies polynomials using the na\"ive $O(n^2)$ algorithm, multiplying each pair of monomials and then summing the results.  However, it is written in C++, and for small polynomials it is superior to \texttt{fftMultiply}, which has the overhead of the top-level \textsc{Macaulay2} language and calls to \textsc{libffi} and \textsc{FFTW}.  But for larger polynomials, \texttt{fftMultiply} can be faster.  For example, in \cref{listing: poly mult compare}, we see that it computed the product of two degree 6000 polynomials over a second more quickly than \textsc{Macaulay2}'s native polynomial multiplication on a system with an AMD Ryzen 5 2600 processor running Ubuntu 22.04.

\begin{listing}[h]
\begin{minted}[frame=lines,framesep=1ex]{macaulay2}
i12 : R = CC[x];

i13 : n = 6000;

i14 : f = (random(CC^1, CC^(n + 1)) * matrix apply(n + 1, i -> {x^i}))_(0,0);

i15 : g = (random(CC^1, CC^(n + 1)) * matrix apply(n + 1, i -> {x^i}))_(0,0);

i16 : elapsedTime f * g;
 -- 16.5132 seconds elapsed

i17 : elapsedTime fftMultiply(f, g);
 -- 15.2714 seconds elapsed
\end{minted}
\caption{Polynomial multiplication running time comparisons}
\label{listing: poly mult compare}
\end{listing}

\section{General Gauss-Markov linear model problem example}
\label{section: glm example}

Statisticians consider the \textit{general Gauss-Markov linear model}
\begin{equation*}
  \vb d = A\boldsymbol\beta + \boldsymbol\varepsilon,
\end{equation*}
with response vector $\vb d\in\RR^n$, $n\times m$ design matrix $A$, parameter vector $\boldsymbol\beta\in\RR^m$, and noise vector $\boldsymbol\varepsilon\in\RR^n$  with mean $\vb 0$ and covariance matrix $\sigma^2W$, where $\sigma^2\in\RR$ is unknown and $W$ is a known $n\times n$ symmetric nonnegative definite matrix of rank $p$.  As shown in \cite{KP}, a \textit{best linear unbiased estimator (BLUE)} of $\boldsymbol\beta$ is the vector $\vb x\in\RR^m$ that, alongside $\vb y\in\RR^p$, minimizes $\|\vb y\|$ (using the usual $L^2$-norm) subject to the constraint $\vb d = A\vb x + B\vb y$, where $B$ is an $n\times p$ matrix of full rank for which $W=BB^T$.

Constrained least squares problems such as this can be solved using the function \texttt{dggglm} from the linear algebra library LAPACK \cite{LAPACK}.  \textsc{Macaulay2} already uses a number of LAPACK functions for its linear algebra computations.  However, \texttt{dggglm} is not one of these.

Using the \texttt{ForeignFunctions} package, it is possible to call this function from \textsc{Macaulay2}.  See \cref{listing: glm functions} for the code.  Note the helper method \texttt{toLAPACK}, which converts \textsc{Macaulay2} matrices and vectors, which are stored using row-major order, into arrays in column-major order as expected by LAPACK.  Note also that since \textsc{Macaulay2} is already linked against LAPACK, it is not necessary to call \texttt{openSharedLibrary} or to include a shared library object as input to \texttt{foreignFunction}.

\begin{listing}[h]
  \begin{minted}[frame=lines,framesep=1ex]{macaulay2}
needsPackage "ForeignFunctions"

toLAPACK = method()
toLAPACK Matrix := A -> (
    T := (numRows A * numColumns A) * double;
    T flatten entries transpose A)
toLAPACK Vector := toLAPACK @@ matrix

dggglm = foreignFunction("dggglm_", void, toList(13:voidstar))

generalLinearModel= method()
generalLinearModel(Matrix, Matrix, Vector) := (A, B, d) -> (
    if numRows A != numRows B
    then error "expected first two arguments to have the same number of rows";
    n := numRows A;
    m := numColumns A;
    p := numColumns B;
    x := getMemory(m * size double);
    y := getMemory(p * size double);
    lwork := n + m + p;
    work := getMemory(lwork * size double);
    i := getMemory int;
    dggglm(address int n, address int m, address int p, toLAPACK A,
	address int n, toLAPACK B, address int n, toLAPACK d, x, y, work,
	address int lwork, i);
    if value(int * i) != 0 then error("call to dggglm failed");
    (vector value (m * double) x, vector value (p * double) y))
  \end{minted}
  \caption{\textsc{Macaulay2} wrapper around \texttt{dggglm} from LAPACK}
  \label{listing: glm functions}
\end{listing}

As an example, we use \textsc{Macaulay2}, via a foreign function call to LAPACK, to reproduce the computation from \cite[\S4]{KP} in \cref{listing: glm example}.

\begin{listing}[H]
  \begin{minted}[frame=lines,framesep=1ex]{macaulay2}
i2 : A = matrix {{1, 2, 3}, {4, 1, 2}, {5, 6, 7}, {3, 4, 6}};

              4       3
o2 : Matrix ZZ  <-- ZZ

i3 : B = matrix {{1, 0, 0, 0}, {2, 3, 0, 0}, {4, 5, 1e-5, 0}, {7, 8, 9, 10}};

                4         4
o3 : Matrix RR    <-- RR
              53        53

i4 : d = vector {1, 2, 3, 4};

       4
o4 : ZZ

i5 : generalLinearModel(A, B, d)

o5 = (|   .3141  |, | .0296627 |)
      | -.334417 |  | .0451036 |
      |  .441691 |  | .0585128 |
                    | .0650142 |

o5 : Sequence
  \end{minted}
  \caption{Example of general Gauss-Markov linear model problem using LAPACK}
  \label{listing: glm example}
\end{listing}
\section{Just-in-time compilation example}
\label{section: jit example}

Another possible use of the \texttt{ForeignFunctions} package is \textit{just-in-time}, or \textit{JIT}, compilation.  At runtime, C code can be compiled into a shared library and called from Macaulay2.  This can result in significant performance increases.

As an example, we consider the na\"ive computation of Fibonacci numbers using their usual recurrence relation definition, i.e., for all nonnegative integers $n$,
\begin{equation*}
  F_n = \begin{cases}
          n &\text{if }n<2\\  
          F_{n-1} + F_{n-2} &\text{otherwise}.
        \end{cases}
      \end{equation*}

      It is a bad idea in general to use this definition to compute Fibonacci numbers due to the fact that its running time is $\Theta(\varphi^n)$, where $\varphi=\frac{1+\sqrt 5}{2}$ is the golden ratio \cite[Section 27.1]{CLRS}.  However, it will be useful to illustrate our point.

      In \cref{listing: fibonacci functions}, we define two functions.  The first, \texttt{fibonacci1}, is a standard top-level Macaulay2 implementation of the above definition.  The second, \texttt{fibonacci2}, uses JIT compilation.  In particular, it writes a C implementation of the definition to a file, calls GCC \cite{GCC} to compile this code into a shared library object, opens this shared library, creates a foreign function, calls it, and finally converts the return value from a C \texttt{int} to a Macaulay2 \texttt{ZZ} object.

\begin{listing}[h]
  \begin{minted}[frame=lines,framesep=1ex]{macaulay2}
needsPackage "ForeignFunctions"

fibonacci1 = n -> if n < 2 then n else fibonacci1(n - 1) + fibonacci1(n - 2)

fibonacci2 = n -> (
    dir := temporaryFileName();
    makeDirectory dir;
    dir | "/libfib.c" << ///int fibonacci2(int n)
{
    if (n < 2)
        return n;
    else
        return fibonacci2(n - 1) + fibonacci2(n - 2);
}
/// << close;
    run("gcc -c -fPIC " | dir | "/libfib.c -o " | dir | "/libfib.o");
    run("gcc -shared " | dir | "/libfib.o -o " | dir | "/libfib.so");
    lib := openSharedLibrary("libfib", FileName => dir | "/libfib.so");
    f = foreignFunction(lib, "fibonacci2", int, int);
    value f n)
  \end{minted}
  \caption{Definitions of Macaulay2 and JIT functions implementing Fibonacci number computation}
  \label{listing: fibonacci functions}
\end{listing}

As can be seen in \cref{listing: fibonacci comparison}, despite the additional overhead of compiling a shared library and making a foreign function call, the JIT implementation is significantly quicker, running over 80 times faster than the native Macaulay2 implementation in this particular test.

\section*{Acknowledgements}

We would like to the thank the anonymous referee for helpful comments, and in particular for suggesting the examples in Sections \ref{section: glm example} and \ref{section: jit example}.

\begin{listing}[H]
  \begin{minted}[frame=lines,framesep=1ex]{macaulay2}
i2 : elapsedTime fibonacci1 35
 -- 10.5189s elapsed

o2 = 9227465

i3 : elapsedTime fibonacci2 35
 -- .129873s elapsed

o3 = 9227465
  \end{minted}
  \caption{Comparison of Macaulay2 and JIT Fibonacci number computation performances}
  \label{listing: fibonacci comparison}
\end{listing}

\bibliography{foreign-functions}{}

\begin{thebibliography}{10}

\bibitem{LAPACK}
E.~Anderson, Z.~Bai, C.~Bischof, L.~S. Blackford, J.~Demmel, J.~Dongarra,
  J.~Du~Croz, A.~Greenbaum, S.~Hammarling, A.~McKenney, et~al.
\newblock {\em LAPACK users' guide}.
\newblock SIAM, 1999.

\bibitem{MPSolve1}
D.~A. Bini and G.~Fiorentino.
\newblock Design, analysis, and implementation of a multiprecision polynomial
  rootfinder.
\newblock {\em Numer. Algorithms}, 23(2-3):127--173, 2000.

\bibitem{MPSolve2}
D.~A. Bini and L.~Robol.
\newblock Solving secular and polynomial equations: A multiprecision algorithm.
\newblock {\em Journal of Computational and Applied Mathematics}, 272:276--292,
  2014.

\bibitem{Normaliz}
W.~Bruns, C.~S. B.~Ichim, and U.~von~der Ohe.
\newblock Normaliz. algorithms for rational cones and affine monoids.
\newblock Available at \url{https://normaliz.uos.de}.

\bibitem{Normaliz-M2}
W.~Bruns and G.~K\"{a}mpf.
\newblock A {M}acaulay2 interface for {N}ormaliz.
\newblock {\em J. Softw. Algebra Geom.}, 2:15--19, 2010.

\bibitem{CT}
J.~W. Cooley and J.~W. Tukey.
\newblock An algorithm for the machine calculation of complex {F}ourier series.
\newblock {\em Math. Comp.}, 19:297--301, 1965.

\bibitem{CLRS}
T.~H. Cormen, C.~E. Leiserson, R.~L. Rivest, and C.~Stein.
\newblock {\em Introduction to algorithms}.
\newblock MIT Press, Cambridge, MA, third edition, 2009.

\bibitem{mpfr}
L.~Fousse, G.~Hanrot, V.~Lef\`{e}vre, P.~P\'{e}lissier, and P.~Zimmermann.
\newblock {MPFR}: A multiple-precision binary floating-point library with
  correct rounding.
\newblock {\em ACM Trans. Math. Softw.}, 33(2):13–es, jun 2007.

\bibitem{FFTW}
M.~Frigo and S.~Johnson.
\newblock The design and implementation of {FFTW3}.
\newblock {\em Proceedings of the IEEE}, 93(2):216--231, 2005.

\bibitem{gmp}
T.~Granlund.
\newblock {GNU} multiple precision arithmetic library.
\newblock Available at \url{https://gmplib.org/}.

\bibitem{M2}
D.~R. Grayson and M.~E. Stillman.
\newblock Macaulay2, a software system for research in algebraic geometry.
\newblock Available at \url{https://macaulay2.com/}.

\bibitem{libffi}
A.~Green.
\newblock libffi, a portable foreign function interface library.
\newblock Available at \url{https://sourceware.org/libffi/}.

\bibitem{KP}
S.~Kourouklis and C.~C. Paige.
\newblock A constrained least squares approach to the general {G}auss-{M}arkov
  linear model.
\newblock {\em J. Amer. Statist. Assoc.}, 76(375):620--625, 1981.

\bibitem{GCC}
R.~Stallman.
\newblock {\em Using GCC: The GNU Compiler Collection: Reference Manual}.
\newblock Free Software Foundation. GNU Press, 2003.

\end{thebibliography}
\bibliographystyle{abbrv}

\end{document}